\documentclass{article}

\usepackage{arxiv}
\usepackage[utf8]{inputenc}
\usepackage[T1]{fontenc}
\usepackage{amsmath,amssymb,amsfonts,mathtools}
\usepackage{amsthm}
\usepackage{microtype}
\usepackage{enumitem}
\usepackage[numbers,sort&compress]{natbib}
\usepackage{url}
\usepackage[hidelinks,
 pdftitle={Moment classes in information cohomology: injectivity into functional coefficients},
 pdfauthor={Vladimir A. Lazovskiy, Victor B. Kazantsev, Sergey V. Stasenko}]{hyperref}

\newtheorem{Theorem}{Theorem}[section]
\newtheorem{Lemma}[Theorem]{Lemma}
\newtheorem{Proposition}[Theorem]{Proposition}
\newtheorem{Corollary}[Theorem]{Corollary}
\theoremstyle{definition}
\newtheorem{Definition}[Theorem]{Definition}
\newtheorem{Remark}[Theorem]{Remark}
\numberwithin{equation}{section}
\newcommand{\R}{\mathbb R}\newcommand{\C}{\mathbb C}\newcommand{\Z}{\mathbb Z}\newcommand{\EE}{\mathbb E}
\newcommand{\Sym}{\operatorname{Sym}}\newcommand{\ev}{\operatorname{ev}}
\newcommand{\LL}{\mathbb L}\newcommand{\Loc}{\mathrm{Loc}}
\newcommand{\Mod}{\mathrm{Mod}}\newcommand{\bA}{\mathbf A}\newcommand{\Ext}{\operatorname{Ext}}
\newcommand{\FF}{\mathcal F}\newcommand{\Dop}{\mathcal L}\newcommand{\Om}{\mathcal U}\newcommand{\cX}{\mathcal Y}

\title{Moment classes in information cohomology:\\ injectivity into functional coefficients}
\author{
  Vladimir A. Lazovskiy$^{1}$
  \And
  Victor B. Kazantsev$^{1,2}$\\
  \texttt{kazantsev@neuro.nnov.ru}
  \And
  Sergey V. Stasenko$^{1,2}$\\
  \texttt{stasenko@neuro.nnov.ru}
}

\begin{document}
\maketitle

\begin{center}
\small
$^{1}$ Moscow Center for Advanced Studies, 20 Kulakova Str., 123592 Moscow, Russia\\
$^{2}$ Lobachevsky State University of Nizhny Novgorod, 23 Gagarin Avenue, 603022 Nizhny Novgorod, Russia
\end{center}

\begin{abstract}
For $n\ge2$ finite random variables, each taking at least two values, we consider the module $\LL^{(m)}$ of degree $m$ polynomials in formal logarithms and its realization by moments of self-information. We prove that the induced map $H^k(S,\LL^{(m)})\to H^k(S,\FF_\Pi)$ into information cohomology with measurable coefficients is injective for all $k\ge0$ and $m\ge1$. Its source is nonzero only for $k=n$ with $m\ge2$ and for $(k,m)=(1,1)$, so these are the cases with content. In particular, an $m$-th moment representing a nonzero formal class in degree $n$ is not a sum of averaged functions of conditional laws and functions of $(n-1)$-variable marginals, even if these functions are arbitrary. The proof uses finite differences, algebraic independence of logarithms of linear forms, and a computation at product laws.
\end{abstract}


\section{Introduction}\label{sec:intro}

The information cohomology $H^k(S,F)=\Ext^k_{\Mod(\bA)}(\R_S,F)$ of Baudot--Bennequin and Vigneaux \cite{BB,V,Vth} attaches groups to a structure $S$ of observables with coefficients in a module $F$, typically a module of functions of the probability law. In degree one the Shannon entropy appears as essentially the only measurable class \cite{BB}, \cite[Theorem~1.3.2]{V}. In higher degrees the cohomology with functional coefficients is much less accessible. In a separate work~\cite{LKS1} a finite-dimensional model is studied: the modules $\LL^{(m)}$ of homogeneous polynomials of degree $m$ in formal symbols $\lambda_T$, on which conditioning acts as on self-informations. For the structure of $n$ jointly measurable variables their cohomology is concentrated in degree $n$ (apart from one line in degree one for $m=1$), and the top group is an explicit quotient $Q_n^{(m)}$ of $\Sym^mV$. Substituting $\lambda_T=-\log P_T$ and averaging gives a morphism $\ev_\Pi$ from $\LL^{(m)}$ to the module $\FF_\Pi$ of measurable functions of the law. The question is whether the formal classes survive in the functional cohomology, where the coboundaries are built from \emph{arbitrary} measurable functions.

\begin{Theorem}\label{thm:A}
Let $n\ge2$, let $X_1,\dots,X_n$ take values in finite sets $\cX_i$ with $|\cX_i|\ge2$, and let $S$ be the Boolean structure of these variables. Then for all $k\ge0$ and $m\ge1$ the map
\[
(\ev_\Pi)_*:\ H^k(S,\LL^{(m)})\longrightarrow H^k(S,\FF_\Pi)
\]
is injective.
\end{Theorem}

By the concentration theorem of \cite{LKS1}, together with the vanishing of $H^n(S,\LL^{(1)})$ proved there, the source is nonzero only for $k=n$ with $m\ge2$ and for $(k,m)=(1,1)$, so these are the cases with content. The case $(1,1)$ concerns the entropy class and follows from known facts; the main point is the top degree $k=n$, where the statement takes the following concrete form. Write $V$ for the span of the $\lambda_T$, $\varnothing\ne T\subseteq [n]=\{1,\dots,n\}$, and $R_m\subseteq\Sym^mV$ for the subspace of local polynomials \eqref{eq:Rm}.

\begin{Corollary}\label{cor:B}
In the setting of Theorem~\ref{thm:A}, let $g\in\Sym^mV$, and consider the equality
\begin{equation}\label{eq:dec}
\EE_P\,g\big((-\log P_T)_T\big)=\sum_{i\in [n]}\sum_{a:\,P_i(a)>0}P_i(a)\,\varphi_i\big(P|_{X_i=a}\big)+\sum_{i\in [n]}\psi_i\big(P_{[n]\setminus i}\big),
\end{equation}
where $\varphi_i$ is a function of a law of $(X_1,\dots,X_n)$ and $\psi_i$ is a function of a law of the other $n-1$ variables. The following are equivalent:
\begin{enumerate}[label=(\roman*)]
\item $g\in R_m$;
\item \eqref{eq:dec} holds for all laws $P$, for some functions $\varphi_i,\psi_i$;
\item \eqref{eq:dec} holds for all positive laws $P$, for some \emph{arbitrary} functions $\varphi_i,\psi_i$.
\end{enumerate}
\end{Corollary}

Thus allowing arbitrary functions, instead of polynomials in self-informations, produces no new relations between moments. We do not claim surjectivity, we do not compute $H^k(S,\FF_\Pi)$ for $k\ge2$, and the statements concern moments, which are linear in $\ev_\Pi$: nothing is claimed about products of entropies, such as the term $H(P)^2$ of the varentropy.

\paragraph{Idea of the proof.} We replace laws by positive unnormalized vectors $x$ on the set of joint outcomes. A function of the form $b_i\varphi$ or $\iota_i\psi$ becomes a sum of functions of ``slices'' of $x$. Constant coefficient differential operators whose symbols vanish on the corresponding subspaces annihilate every real-analytic sum of such functions, although the summands are arbitrary; this is proved with finite differences. Applied to the moment function, such an operator yields a polynomial identity between logarithms of linear forms, which by monodromy must hold identically. The resulting linear conditions on $g$ are then evaluated at one product law, where they are strong enough to span the annihilator of $R_m$.

\paragraph{Organization.} Section~\ref{sec:recall} recalls the facts from \cite{LKS1} that are used; none of the results of the present paper is used in \cite{LKS1}. Section~\ref{ssec:Fcomp} introduces the functional modules. Section~\ref{app:inj} proves the injectivity in the top degree, Section~\ref{sec:cons} deduces Theorem~\ref{thm:A} and discusses degree one, and Section~\ref{sec:open} lists open questions.

\section{Notation and the facts used}\label{sec:recall}

Throughout, $n\ge2$, $[n]=\{1,\dots,n\}$, and $S$ is the Boolean information structure of $n$ variables: its objects are the subsets $y\subseteq [n]$ (the joint measurement of $X_i$, $i\in y$), with an arrow $y\to y'$ if $y'\subseteq y$ \cite[Ex.~1.2.2]{V}. The coefficient category $\Mod(\bA)$ consists of modules over the sheaf of monoid algebras of $S$ \cite[\S3.1]{V}; explicitly, a module $G$ over $S$ consists of spaces $G(y)$ with commuting idempotent operators $b_x$, $x\subseteq y$, such that $b_\varnothing=\mathrm{id}$ and $b_xb_{x'}=b_{x\cup x'}$, and of maps $\iota:G(y')\to G(y)$, $y'\subseteq y$, compatible with the $b_x$ and functorial. We write $b_i=b_{\{i\}}$, $e_x=1-b_x$, and $\iota_i:G([n]\setminus i)\to G([n])$.

\subsection{The formal modules} Let $V_y$ be the space with basis $\lambda_T$, $\varnothing\ne T\subseteq y$, put $\lambda_\varnothing=0$, $V=V_{[n]}$, and let $b_x\lambda_T=\lambda_{T\cup x}-\lambda_x$. The module $\LL^{(m)}$ is $\LL^{(m)}(y)=\Sym^mV_y$, with $b_x$ acting as ring endomorphisms of $\Sym V_y$ and with the natural embeddings {\cite[Definition~14]{LKS1}}. Put
\begin{equation}\label{eq:Rm}
R_m=\sum_{i\in [n]}b_i\,\Sym^mV+\sum_{i\in [n]}\Sym^mV_{[n]\setminus i},\qquad Q_n^{(m)}=\Sym^mV/R_m .
\end{equation}

\subsection{Cohomology} We use the following facts.

\begin{Proposition}[reduced complex, {\cite[Corollary~2]{LKS1}}]\label{cor:reduced}
For every module $G$ the groups $H^k(S,G)$ are computed, naturally in $G$, by the complex with terms $\tilde C^k=\bigoplus_{A\subseteq [n],\ |A|=k}\tilde e_AG(A)$ and differential $(\delta f)_A=\sum_j\pm\tilde e_A\,\iota f_{A\setminus a_j}$, where $\tilde e_A=\prod_{a\in A}(1-b_a)$ (not to be confused with $1-b_A$).
\end{Proposition}

\begin{Proposition}[top degree, {\cite[Proposition~10]{LKS1}}]\label{inj:propA}
For every module $G$ there is an isomorphism, natural in $G$,
\[
H^n(S,G)\cong G([n])\Big/\Big(\sum_ib_iG([n])+\sum_i\iota_iG([n]\setminus i)\Big).
\]
In particular $H^n(S,\LL^{(m)})\cong Q_n^{(m)}$ {\cite[Corollary~5]{LKS1}}.
\end{Proposition}

\begin{Theorem}[concentration, {\cite[Theorem~4]{LKS1}}]\label{thm:conc}
For $m\ge1$ one has $H^k(S,\LL^{(m)})=0$ for $k\ne n$, except that $H^1(S,\LL^{(1)})\cong\R$. \end{Theorem}

\noindent Moreover $Q_n^{(1)}=0$, so that $H^n(S,\LL^{(1)})=0$ {\cite[Corollary~5]{LKS1}}.

\subsection{The cube model and the dual of the top group} Let $h_C$, $C\subseteq [n]$, be the function on $\{0,1\}^n$ equal to $1$ on the points with support in $C$ and $0$ elsewhere. The map
\begin{equation}\label{b:cube}
\lambda_T\mapsto h_{[n]}-h_{[n]\setminus T}
\end{equation}
identifies $V$ with the space of functions vanishing at $0$; under it $b_i$ is multiplication by $1-s_i$, where $s_i$ is the $i$-th coordinate on $\{0,1\}^n$. Let $\nu_C$ (denoted $y_C$ in \cite{LKS1}) be the coordinates of $V$ with respect to the functions $h_C$; they satisfy $\sum_C\nu_C=0$ and generate the algebra of polynomials on $V$. The dual of $\Sym^mV$ is the space of homogeneous polynomials of degree $m$ on $V$, and
\begin{equation}\label{b:ann}
R_m^{\perp}=\{\Psi:\ \Psi|_{b_iV}=0,\ \Psi|_{V_{[n]\setminus i}}=0\ \ \forall i\}.
\end{equation}
A multiset $\chi$ of $m$ subsets of $[n]$ is \emph{essential} if every $i\in [n]$ belongs to some but not to all members of $\chi$; put $\nu^\chi=\prod_{C\in\chi}\nu_C$.

\begin{Proposition}[dual description, {\cite[Proposition~11]{LKS1}}]\label{b:dual}
For $m\ge2$ the space $R_m^\perp$ is spanned by the restrictions to $V$ of the monomials $\nu^\chi$ with $\chi$ essential, $|\chi|=m$.
\end{Proposition}

\section{Modules of functions of the law}\label{ssec:Fcomp}

\begin{Definition}[the module $\FF$ of functions of the law]\label{inj:defF}
For $y\subseteq [n]$ let $\mathrm{Box}(y)$ be the set of laws $P$ on $\cX_y=\prod_{i\in y}\cX_i$ whose support is a box $\prod_{i\in y}\cX'_i$ ($\varnothing\ne \cX'_i\subseteq \cX_i$) and which are positive on it. This set contains all positive laws, and it is closed under marginals and under conditioning $P\mapsto P|_{X_x=a}$ ($P_x(a)>0$). Put $\FF(y)=\R^{\mathrm{Box}(y)}$ (all real functions of $P\in\mathrm{Box}(y)$), let $\iota:\FF(y')\to\FF(y)$ for $y'\subseteq y$ be the pullback of the marginal, $(\iota\psi)(P)=\psi(P_{y'})$, and let $x\subseteq y$ act by
\[
(b_x\varphi)(P)=\sum_{a:\,P_x(a)>0}P_x(a)\,\varphi(P|_{X_x=a}).
\]
The moment map $\ev:\LL^{(m)}(y)\to\FF(y)$ is $\ev(g)(P)=\EE_P\,g\big((\imath_T)_T\big)$, $\imath_T=-\log P_T(X_T)$.
\end{Definition}
\begin{Lemma}\label{inj:mod}
(a) $\FF$ is a module over $S$: $b_xb_{x'}=b_{x\cup x'}$, $b_\varnothing=\mathrm{id}$ and $b_x\iota=\iota b_x$. (b) $\ev:\LL^{(m)}\to\FF$ is a morphism of modules.
\end{Lemma}
\begin{proof}
(a) In $(b_xb_{x'}\varphi)(P)=\sum_aP_x(a)\sum_{a'}(P|_a)_{x'}(a')\varphi(P|_{a,a'})$ only the $a'$ agreeing with $a$ on $x\cap x'$ contribute, and $P_x(a)(P|_a)_{x'}(a')=P_{x\cup x'}(a'')$ for the glued value $a''$. Commutation with $\iota$: $(P|_a)_{y'}=P_{y'}|_a$ for $x\subseteq y'$. (b) On the support of $P|_{X_x=a}$ we have $\imath_C(P|_a)=\imath_{C\cup x}(P)-\imath_x(P)$, i.e.\ the value of $b_x\lambda_C$; hence $\ev(b_xg)(P)=\sum_aP_x(a)\,\ev(g)(P|_a)$. Compatibility with $\iota$: $\imath_T$ for $T\subseteq y'$ depends only on $P_{y'}$.
\end{proof}
Positive laws on the full alphabets do not form a module: $b_x$ maps $\FF(y)$ to itself, so $\varphi$ must be evaluated at the conditional law $P|_{X_x=a}$ regarded as a law on $\cX_y$, which is concentrated on $\{X_x=a\}$; its support is a box, but it is not positive. This is why laws with box supports are used. The module $\FF$ is auxiliary; the module of~\cite{V} is the following one.

\begin{Definition}[the module $\FF_\Pi$ of measurable functions on all laws]\label{inj:defPi}
For $y\subseteq [n]$ let $\Pi(y)$ be the simplex of all laws on $\cX_y$ (zeros are allowed) and let $\FF_\Pi(y)$ be the space of Borel measurable functions $\Pi(y)\to\R$; the map $\iota$ is the pullback of the marginal, and
\[
(b_x\varphi)(P)=\sum_{a:\ P_x(a)>0}P_x(a)\,\varphi(P|_{X_x=a}).
\]
This is the module $\mathcal F_1(Q)$ of~\cite{V} (exponent $\alpha=1$) for the functor $Q$ of all laws; that it is a module is \cite[Prop.~4.2.1, 4.2.2]{V}. For $g\in\LL^{(m)}(y)$ put
\[
\ev_\Pi(g)(P)=\sum_{u:\ P(u)>0}P(u)\,g\big((\imath_T(u))_T\big),\qquad \imath_T(u)=-\log P_T(u_T),\quad u\in \cX_y.
\]
\end{Definition}
\begin{Lemma}\label{inj:Pi}
(a) $\ev_\Pi(g)$ is continuous on $\Pi(y)$, in particular measurable. (b) $\ev_\Pi:\LL^{(m)}\to\FF_\Pi$ is a morphism of modules. (c) The restriction $\mathrm{res}:\FF_\Pi\to\FF$ to $\mathrm{Box}(y)\subseteq\Pi(y)$ is a morphism of modules and $\mathrm{res}\circ\ev_\Pi=\ev$.
\end{Lemma}
\begin{proof}
(a) If $P(u)>0$, then $P_T(u_T)\ge P(u)>0$ for all $T\subseteq y$, hence $0\le\imath_T(u)\le-\log P(u)$. Since $g$ is a homogeneous polynomial of degree $m$ in the $\imath_T$, we get $|P(u)\,g|\le C\,P(u)\,|\log P(u)|^m\to0$ as $P(u)\to0$. So each summand $P\mapsto P(u)\,g(\dots)$ (put equal to $0$ where $P(u)=0$) is continuous, and the sum is finite. (b) The computation of Lemma~\ref{inj:mod}(b) applies verbatim: all sums run over points (respectively values $a$) of positive probability, and the terms with $P_x(a)=0$ are absent on both sides. (c) $\mathrm{Box}(y)$ is closed under marginals and conditioning, and the formulas for $\iota$ and $b_x$ are the same; on $\mathrm{Box}(y)$ the formula for $\ev_\Pi$ is that for $\ev$.
\end{proof}

\section{Injectivity in the top degree}\label{app:inj}

\begin{Proposition}[injectivity into the cohomology with functional coefficients]\label{inj:A}
Let $n\ge2$ and let $\cX_1,\dots,\cX_n$ be finite sets with $|\cX_i|\ge2$. For every $m\ge1$ the map
\[
\ev_*:\ Q_n^{(m)}=H^n(S,\LL^{(m)})\longrightarrow H^n(S,\FF)
\]
induced by the morphism of modules $\ev:\LL^{(m)}\to\FF$ is injective.
\end{Proposition}
By Proposition~\ref{inj:propA} and its naturality, the map $\ev_*$ is $g\mapsto[\ev(g)]$, and its kernel is
\begin{equation}\label{inj:ker}
\{g\in\Sym^mV:\ \ev(g)\in\textstyle\sum_ib_i\FF([n])+\sum_i\iota_i\FF([n]\setminus i)\}\big/R_m .
\end{equation}
It is convenient to work over $\C$; all the maps are defined over $\R$, so injectivity of the complexified map implies injectivity of the real one.

\paragraph{Plan.} Let $g\in\Sym^mV$ be a class whose image in $H^n(S,\FF)$ vanishes. By Proposition~\ref{inj:propA} this means that the function $P\mapsto\ev_P(g)$ is a sum of ``local'' functions: functions $b_i\varphi$ (an average over $X_i$ of a function of the conditional law $P|_{X_i=a}$) and functions $\iota_i\psi$ (a function of the marginal law on $[n]\setminus i$), where $\varphi,\psi$ are \emph{arbitrary} functions. We have to show that $g\in R_m$. The argument goes as follows.
\begin{enumerate}
\item[(1)] \emph{Unnormalized laws.} Replace a law $P$ by a positive vector $x=(x_c)$ indexed by the cells $c=(c_1,\dots,c_n)$ (the possible joint outcomes), $P=x/N$, $N=\sum_cx_c$. A local function becomes a sum of functions each of which depends only on a ``slice'' $B_jx$ of $x$: the restriction of $x$ to a hyperplane $\{c_i=a\}$, or the sum of $x$ over the coordinate $c_i$ (\S\ref{inj:s2}).
\item[(2)] \emph{Operators that kill slices.} A constant coefficient differential operator $D(\partial)$ whose symbol $D$ vanishes on the ``slice subspaces'' $W_j$ kills a function of $B_jx$ by the chain rule. The functions of the sum are arbitrary and not differentiable, so we use finite differences instead; the result is that such $D(\partial)$ kills every \emph{analytic} sum of slice functions (\S\ref{inj:s3}).
\item[(3)] \emph{Explicit form of the function.} The function $\widehat G_g(x)=N\,\ev(g)(x/N)$ has the explicit form $\sum_cx_c\,g(y_{c,\bullet}(x))$, where $y_{c,T}(x)$ are logarithms of linear forms in $x$. Logarithms of pairwise non-proportional linear forms are algebraically independent (monodromy, \S\ref{inj:s4}); therefore $D(\partial)\widehat G_g=0$ forces the polynomial coefficients to vanish. The outcome is a family of linear conditions $\langle\Psi,g\rangle=0$, $\Psi\in Z_m(x)$, where $Z_m(x)$ is an explicit space of polynomials depending on the point $x$.
\item[(4)] \emph{Criterion.} By duality $R_m=I_m^{\perp}$, where $I_m$ is the explicit space spanned by the ``essential'' monomials. So it is enough to find \emph{one} point $x$ with $I_m\subseteq Z_m(x)$: then $g\perp I_m$, that is $g\in R_m$ (\S\ref{inj:s4}).
\item[(5)] \emph{Product points.} For $x_c=\prod_ip_i(c_i)$ and operators $D$ of tensor form the conditions can be computed in closed form; the computation shows that they produce all essential monomials (\S\ref{inj:s5}).
\end{enumerate}

\paragraph{Notation used in this section.}
\begin{itemize}[leftmargin=2em]
\item $\cX=\prod_i\cX_i$ is the set of \emph{cells} $c=(c_1,\dots,c_n)$; $\Om=(0,\infty)^{\cX}$ is the set of positive vectors $x=(x_c)$ (unnormalized laws); $N=\sum_cx_c$; for $\varnothing\ne T\subseteq [n]$ and $a\in \cX_T=\prod_{i\in T}\cX_i$ the \emph{marginal sums} are $q_{T,a}(x)=\sum_{c:\,c|_T=a}x_c$ (so $q_\varnothing=N$).
\item $y_{c,T}(x)=\log N-\log q_{T,c|_T}(x)$ is the self-information of the cell $c$ with respect to $T$ for the law $P=x/N$. The same expression with an arbitrary $a\in \cX_T$ in place of $c|_T$ is denoted $z_{T,a}(x)$.
\item $\Sym^mV=\C[\lambda_T]_m$. When $g\in\Sym^mV$ is evaluated at a point we write $g(y)$, $y=(y_T)_T$: the variable $y_T$ stands for $\lambda_T$. The polynomial algebra $\C[V]=\C[\mu_T]$ is the algebra of polynomials in the coordinates $\mu_T$ dual to $\lambda_T$; it pairs with $\Sym^mV$ by $\langle\Psi,g\rangle=\Psi(\partial_\lambda)g$, so $\mu_T$ acts as $\partial/\partial y_T$.
\item $\nu_C$ ($C\subseteq [n]$) are the coordinates of the cube model \eqref{b:cube}; they are linear forms in $\mu$ with $\sum_C\nu_C=0$ (see the beginning of \S\ref{inj:s4}). Multisets of subsets of $[n]$ are denoted by $\chi$, and $\nu^\chi=\prod_{C\in\chi}\nu_C$.
\end{itemize}

\subsection{Reduction to local functions}\label{inj:s2}
We use the notation $\cX$, $\Om$, $N$, $q_{T,a}$, $y_{c,T}$ listed above: for the law $P=x/N$ we have $-\log P_T(c|_T)=y_{c,T}(x)$. For $g\in\Sym^mV=\C[\lambda_T]_m$ put
\[
\widehat G_g(x)=N\,\ev(g)(x/N)=\sum_{c\in\cX}x_c\,g\big(y_{c,\bullet}(x)\big),
\]
a real-analytic function on $\Om$. Let $B_{i,a}:\R^{\cX}\to\R^{\{c:\,c_i=a\}}$ be the restriction to the slice $\{c_i=a\}$ and $B_i:\R^{\cX}\to\R^{\prod_{j\ne i}\cX_j}$ the summation over $X_i$; these are $\sum_i(|\cX_i|+1)$ coordinate maps, which we denote $B_j$ ($j\in J$). Let $\Loc^\Om$ be the space of functions $\sum_jf_j\circ B_j$ on $\Om$ with arbitrary $f_j$.
\begin{Lemma}\label{inj:red}
If $\ev(g)\in\sum_ib_i\FF([n])+\sum_i\iota_i\FF([n]\setminus i)$, then $\widehat G_g\in\Loc^\Om$.
\end{Lemma}
\begin{proof}
If $\varphi=b_i\varphi'$, then $\varphi(P)=\sum_aP_i(a)\varphi'(P|_{X_i=a})$ and $P|_{X_i=a}$ is determined by the slice of $P$ over $a$; so $N\varphi(x/N)=\sum_aq_{i,a}\,\varphi'_a\big((x_c)_{c_i=a}/q_{i,a}\big)$ is a sum of functions of the slices $B_{i,a}x$. If $\varphi=\iota_i\psi$, then $N\varphi(x/N)=N\psi(B_ix/N)$ is a function of $B_ix$ ($N$ is the sum of the coordinates of $B_ix$). The restriction to positive laws only discards values.
\end{proof}
By \eqref{inj:ker} and Lemma~\ref{inj:red}, Proposition~\ref{inj:A} follows from
\[(\ast_m)\qquad g\in\Sym^mV,\ \ \widehat G_g\in\Loc^\Om\ \Longrightarrow\ g\in R_m.\]

\subsection{Differential operators that kill local functions}\label{inj:s3}
Let $W_{i,a}\subseteq\C^{\cX}$ be the space of vectors supported on the slice $\{c_i=a\}$, and $W_i\subseteq\C^{\cX}$ the space of vectors that depend only on $c|_{[n]\setminus i}$ (are constant along $c_i$). These are the images of the transposed maps $B_{i,a}^T$ and $B_i^T$; with the index set $J$ of \S\ref{inj:s2} we write $W_j=\operatorname{im}B_j^T$, $j\in J$. Let $\mathfrak S_d$ be the space of homogeneous polynomials $D(\xi)$ of degree $d$ on $\C^{\cX}$ vanishing on $\bigcup_jW_j$, and let $D(\partial)$ be the constant coefficient operator obtained by substituting $\xi_c\mapsto\partial/\partial x_c$.

If $f$ were differentiable, then for $D\in\mathfrak S_d$ the chain rule would give $D(\partial)(f\circ B_j)=\big((D\circ B_j^T)(\partial)f\big)\circ B_j=0$, because $D$ vanishes on $W_j=\operatorname{im}B_j^T$. But the functions $f_j$ in $\Loc^\Om$ are \emph{arbitrary}: only their sum is assumed to be analytic. This is why finite differences are used below.
\begin{Proposition}\label{inj:kill}
If $h\in\Loc^\Om$ is real-analytic and $D\in\mathfrak S_d$, then $D(\partial)h=0$.
\end{Proposition}
\begin{proof}
\emph{Step 1: difference operators kill slice functions.} Let $\Lambda=\Z^{\cX}$ and let $\C[\Lambda]$ be its group algebra (Laurent polynomials in $t_c$, $c\in\cX$; $t^\nu=\prod_ct_c^{\nu_c}$). To $\mathsf K=\sum_\nu w_\nu t^\nu\in\C[\Lambda]$ (a finite sum) and $\eta>0$ associate the operator $\mathsf K_\eta f(x)=\sum_\nu w_\nu f(x+\eta\nu)$; it is defined at $x\in\Om$ for all small $\eta$, since $\Om$ is open. Let $J_j=\ker\big(\C[\Lambda]\to\C[B_j\Lambda]\big)$, $t^\nu\mapsto t^{B_j\nu}$. If $\mathsf K\in J_j$, grouping the terms by $\kappa=B_j\nu$ gives
\[
\mathsf K_\eta(f\circ B_j)(x)=\sum_\kappa\Big(\sum_{\nu:\,B_j\nu=\kappa}w_\nu\Big)f(B_jx+\eta\kappa)=0,
\]
since the inner sums are the coefficients of the image of $\mathsf K$ in $\C[B_j\Lambda]$, which is $0$. Hence $\mathsf K_\eta$ kills $\Loc^\Om$ for every $\mathsf K\in\mathfrak a:=\bigcap_jJ_j$.

\emph{Step 2: from differences to $D(\partial)$.} \emph{Claim:} for every $D\in\mathfrak S_d$ there is $\mathsf K\in\mathfrak a$ with $\widehat{\mathsf K}(\xi):=\sum_\nu w_\nu e^{\nu\cdot\xi}=D(\xi)+O(\xi^{d+1})$. Granting the claim, Taylor's formula for the analytic function $h$ gives, for small $\eta$,
\[
\mathsf K_\eta h(x)=\sum_\nu w_\nu h(x+\eta\nu)=\sum_{k\ge0}\eta^k\sum_\nu w_\nu\frac{(\nu\cdot\partial)^k}{k!}h(x)=\sum_{k\ge0}\eta^k\,\widehat{\mathsf K}_k(\partial)h(x),
\]
where $\widehat{\mathsf K}_k$ is the homogeneous part of degree $k$ of $\widehat{\mathsf K}$. By the claim $\widehat{\mathsf K}_k=0$ for $k<d$ and $\widehat{\mathsf K}_d=D$. By Step 1 the left side is $0$, so $0=\eta^dD(\partial)h(x)+O(\eta^{d+1})$, hence $D(\partial)h(x)=0$.

\emph{Step 3: proof of the claim.} Geometrically: $J_j$ is the ideal of a subtorus $T_j$ of the torus $\operatorname{Spec}\C[\Lambda]=(\C^*)^{\cX}$ whose tangent space at the identity is $W_j$; $\mathfrak a$ is the ideal of $\bigcup_jT_j$. The ideal $\mathfrak S=\bigoplus_d\mathfrak S_d$ is the ideal of $\bigcup_jW_j$. The claim says that the ``tangent cone'' of $\bigcup_jT_j$ at the identity is $\bigcup_jW_j$; this is checked in the formal neighbourhood of the identity, where $t_c=e^{\xi_c}$ turns subtori into linear subspaces.

Put $\Lambda_j=\ker B_j\cap\Lambda$. Since $B_j$ has integer entries, $\Lambda_j\otimes\C=\ker B_j=W_j^{\perp}$ (the annihilator with respect to the pairing $\nu\cdot\xi$), and $\Lambda/\Lambda_j\cong B_j\Lambda$ is torsion free; hence $J_j=(t^\nu-1:\nu\in\Lambda_j)$ is the ideal of a subtorus. Let $\mathfrak m=(t_c-1)\subset\C[\Lambda]$ be the maximal ideal of the identity, and $\widehat R=\C[[\xi]]$ the $\mathfrak m$-adic completion, $t_c=e^{\xi_c}$. The associated graded ring $\operatorname{gr}_{\mathfrak m}\C[\Lambda]$ is the polynomial ring $\C[\xi]$ ($t_c-1\mapsto\xi_c$). For an ideal $\mathfrak b$ of $\C[\Lambda]$, of $\widehat R$ or of $\C[\xi]$, $\operatorname{in}(\mathfrak b)\subseteq\C[\xi]$ denotes the ideal spanned by the lowest degree forms of its elements (with respect to the $\mathfrak m$-adic filtration, resp.\ the degree). Let $\mathfrak l_j\subset\C[\xi]$ be the ideal generated by the linear forms vanishing on $W_j$; then $\mathfrak S=\bigcap_j\mathfrak l_j$.
\begin{enumerate}
\item[(i)] Since $e^{\nu\cdot\xi}-1=(\nu\cdot\xi)\cdot(\text{unit of }\widehat R)$ and $\{\nu\cdot\xi:\nu\in\Lambda_j\}$ spans the linear forms vanishing on $W_j$, we get $J_j\widehat R=\mathfrak l_j\widehat R$.
\item[(ii)] $\widehat R$ is flat over $\C[\Lambda]$ and over $\C[\xi]$, and finite intersections of ideals commute with flat base change. Hence $\mathfrak a\widehat R=\bigcap_jJ_j\widehat R=\bigcap_j\mathfrak l_j\widehat R=\mathfrak S\widehat R$. (This is the point where the passage from the geometry of the subtori to the linear algebra of the tangent spaces takes place.)
\item[(iii)] For an ideal $\mathfrak b\subseteq\C[\Lambda]$ one has $\C[\Lambda]/(\mathfrak b+\mathfrak m^k)=\widehat R/(\mathfrak b\widehat R+\widehat{\mathfrak m}^k)$ for all $k$, so $\operatorname{in}(\mathfrak b)=\operatorname{in}(\mathfrak b\widehat R)$; in particular $\operatorname{in}(\mathfrak a)=\operatorname{in}(\mathfrak a\widehat R)$.
\item[(iv)] $\operatorname{in}(\mathfrak S\widehat R)=\mathfrak S$, because $\mathfrak S$ is a homogeneous ideal of $\C[\xi]$.
\end{enumerate}
Together, $\operatorname{in}(\mathfrak a)=\mathfrak S$. Thus every $D\in\mathfrak S_d$ is the lowest form of some $\mathsf K\in\mathfrak a$, so $\mathsf K\in\mathfrak m^d$ and $\mathsf K=D(t-1)+(\text{terms of order}>d\text{ in }t-1)$. Since $t_c-1=\xi_c+O(\xi^2)$, we get $\widehat{\mathsf K}(\xi)=\mathsf K(e^{\xi})=D(\xi)+O(\xi^{d+1})$, which is the claim.
\end{proof}

\subsection{Constraint polynomials and the criterion}\label{inj:s4}
\begin{Lemma}\label{inj:log}
Let $\ell_1,\dots,\ell_k$ be pairwise non-proportional nonzero linear forms on $\C^{\cX}$ which are positive on $\Om$. Then $\log\ell_1,\dots,\log\ell_k$ (real on $\Om$) are algebraically independent over $\C(x)$.
\end{Lemma}
\begin{proof}
Let $P(x;L_1,\dots,L_k)\ne0$ be a polynomial of minimal degree in $L$ with coefficients in $\C(x)$ such that $P(x;\log\ell_1(x),\dots,\log\ell_k(x))=0$ on $\Om$; multiplying by a common denominator, we may assume that the coefficients are polynomials in $x$. Let $\Sigma\subset\C^{\cX}$ be the union of the hyperplanes $\{\ell_j=0\}$, and $U=\C^{\cX}\setminus\Sigma$; $U$ is connected (the complement of a proper algebraic subset of $\C^{\cX}$), and $\Om\subset U$. The identity $P(x;\log\ell_1(x),\dots,\log\ell_k(x))=0$ continues analytically along any path in $U$, the logarithms being continued along the path. Fix $i$ and take a path in $U$ from a base point of $\Om$ to a point near a generic point of $\{\ell_i=0\}$, i.e.\ a point of $\{\ell_i=0\}$ lying on no other hyperplane $\{\ell_j=0\}$ (it exists since the forms are non-proportional); go once around $\{\ell_i=0\}$ by a small loop there and return along the same path. Along this loop $\log\ell_i\mapsto\log\ell_i+2\pi\sqrt{-1}$, all the other $\log\ell_j$ return to themselves, and the coefficients (polynomials) are unchanged. Thus $Q(\tau)=P(x;\dots,\log\ell_i+\tau,\dots)$ vanishes for $\tau\in2\pi\sqrt{-1}\Z$, hence $Q\equiv0$ as a polynomial in $\tau$; in particular $\partial P/\partial L_i$ vanishes at $L=\log\ell$, and by minimality of the degree $\partial P/\partial L_i=0$. As $i$ is arbitrary, $P$ does not depend on $L$, so $P\in\C(x)$ and $P=0$, a contradiction.
\end{proof}
\begin{Corollary}\label{inj:z}
The functions $z_{T,a}(x)=\log N-\log q_{T,a}(x)$ ($\varnothing\ne T\subseteq [n]$, $a\in \cX_T$) are algebraically independent over $\C(x)$.
\end{Corollary}
\begin{proof}
The forms $q_{T,a}$ and $N$ are pairwise non-proportional because their supports are different; this uses that every $|\cX_i|\ge2$ (for $T\ne\varnothing$ the support of $q_{T,a}$ is a proper nonempty subset of $\cX$, and distinct $(T,a)$ give distinct supports). Apply Lemma~\ref{inj:log} to the forms $N$ and $q_{T,a}$, and note that $z_{T,a}=\log N-\log q_{T,a}$ is an invertible triangular change of the generators.
\end{proof}

\medskip\noindent\emph{The derivative formula.} We compute $D(\partial)\widehat G_g$, $D\in\mathfrak S$, in terms of $g$.

\emph{Coordinates on $V$.} Let $\mu_T$ ($\varnothing\ne T\subseteq [n]$) be the coordinates on $V$ dual to the basis $\lambda_T$, and $\nu_C$ ($C\subseteq [n]$) the coordinates with respect to the basis $h_C$ of the cube model \eqref{b:cube}. Under the identification \eqref{b:cube} the vector $\sum_T\mu_T\lambda_T$ becomes the function $(\sum_T\mu_T)h_{[n]}-\sum_T\mu_Th_{[n]\setminus T}$, so
\[
\nu_C=-\mu_{[n]\setminus C}\ (C\ne [n]),\qquad\nu_{[n]}=\textstyle\sum_T\mu_T,\qquad\textstyle\sum_C\nu_C=0 .
\]
Hence every polynomial in the $\nu_C$ is an element of $\C[V]$, and the only relation among the $\nu_C$ is $\sum_C\nu_C=0$.

\emph{Step 1: a generating function.} For $c\in\cX$ and $C\subseteq [n]$ let $\ell_{c,C}(u)=\sum_{c':\,c'|_{[n]\setminus C}=c|_{[n]\setminus C}}u_{c'}$ (a linear form in $u\in\C^{\cX}$) and $q_{c,C}=\ell_{c,C}(x)$. Then $q_{c,\varnothing}=x_c$, $q_{c,[n]}=N$ and $q_{c,C}=q_{[n]\setminus C,\,c|_{[n]\setminus C}}$ for $C\ne [n]$; thus $y_{c,T}(x)=\log q_{c,[n]}-\log q_{c,[n]\setminus T}$. For $s=(s_T)$ put $\mu_T=s_T$ (so that $\nu_C=\nu_C(s)$ are linear forms in $s$) and
\[
F_c(s;x):=x_c\exp\Big(\sum_Ts_Ty_{c,T}(x)\Big)=x_c\prod_{C\subseteq [n]}q_{c,C}(x)^{\nu_C(s)} .
\]
Since $g(\partial_s)e^{s\cdot y}\big|_{s=0}=g(y)$ for a polynomial $g$,
\[
\widehat G_g(x)=\sum_cx_c\,g\big(y_{c,\bullet}(x)\big)=g(\partial_s)\Big[\sum_cF_c(s;x)\Big]_{s=0}.
\]
\emph{Step 2: derivatives in $x$.} Since $\ell_{c,C}(x+u)=q_{c,C}(x)+\ell_{c,C}(u)$,
\[
F_c(s;x+u)=F_c(s;x)\,(1+u_c/x_c)\prod_C\big(1+\ell_{c,C}(u)/q_{c,C}\big)^{\nu_C}.
\]
Apply $D(\partial_u)$ at $u=0$:
\[
D(\partial_x)F_c(s;x)=\Theta_{c,D}(\nu(s);x)\,\frac{F_c(s;x)}{x_c},
\]
\[
\Theta_{c,D}(\nu;x):=D(\partial_u)\Big[(x_c+u_c)\prod_C\big(\ell_{c,C}(x+u)/q_{c,C}(x)\big)^{\nu_C}\Big]_{u=0}.
\]
Here $\Theta_{c,D}$ is a polynomial in the $\nu_C$ (expand the powers in the binomial series; since $D(\partial_u)$ has finite order only finitely many terms contribute) whose coefficients are rational functions of $x$. Note that $F_c/x_c=e^{s\cdot y_{c,\bullet}(x)}$.

\emph{Step 3: back to $g$.} Write $\Theta_{c,D}(\nu(s);x)=\Upsilon(s)$, a polynomial in $s$. For every polynomial $g$ one has $g(\partial_s)\big[\Upsilon(s)e^{s\cdot y}\big]_{s=0}=(\Upsilon(\partial_y)g)(y)$ (both sides are linear in $g$ and in $\Upsilon$; for monomials this is the Leibniz rule: $\partial_s^\alpha[s^\beta e^{s\cdot y}]_{s=0}=\frac{\alpha!}{(\alpha-\beta)!}y^{\alpha-\beta}=\partial_y^\beta y^\alpha$). Therefore
\begin{equation}\label{inj:main}
D(\partial)\widehat G_g(x)=\sum_c\big(\Theta_{c,D}(\partial_y;x)\,g\big)\big(y_{c,\bullet}(x)\big),
\end{equation}
where in $\Theta_{c,D}(\nu;x)$ the $\nu_C$ are replaced by the above linear forms in $\mu_T=\partial/\partial y_T$.

For $x\in\Om$ let $Z_m(x)\subseteq\C[V]_m$ be the span of the degree $m$ homogeneous components of the polynomials $\sum_c\Theta_{c,D}(\,\cdot\,;x)$, $D\in\mathfrak S=\bigoplus_d\mathfrak S_d$. Recall that $\C[V]_m$ pairs with $\Sym^mV$ by $\langle\Psi,g\rangle=\Psi(\partial_\lambda)g$.
\begin{Proposition}\label{inj:C}
If $g\in\Sym^mV$ and $\widehat G_g\in\Loc^\Om$, then $\langle\Psi,g\rangle=0$ for all $\Psi\in Z_m(x)$ and all $x\in\Om$.
\end{Proposition}
\begin{proof}
The function $\widehat G_g$ is real-analytic and lies in $\Loc^\Om$, so by Proposition~\ref{inj:kill} the left side of \eqref{inj:main} is $0$. Put $P_c(x;y):=\Theta_{c,D}(\partial_y;x)\,g$: a polynomial in $y$ whose coefficients are rational functions of $x$. Then \eqref{inj:main} says $\sum_cP_c\big(x;y_{c,\bullet}(x)\big)=0$ on $\Om$. Since $y_{c,T}(x)=z_{T,c|_T}(x)$, this reads $P(x;z(x))=0$ for the polynomial $P(x;z)=\sum_cP_c\big(x;(z_{T,c|_T})_T\big)\in\C(x)[z_{T,a}]$. By Corollary~\ref{inj:z}, $P=0$ in $\C(x)[z]$; in particular its constant term in $z$ is $0$. The constant term of $P_c(x;(z_{T,c|_T})_T)$ is $P_c(x;0)$. Write $\Theta_{c,D}=\sum_\beta\theta_\beta(x)\mu^\beta$ as a polynomial in $\mu$. As $g$ is homogeneous of degree $m$, only the terms with $|\beta|=m$ contribute to $P_c(x;0)=\sum_\beta\theta_\beta\,\partial^\beta g(0)$, and this is $\langle(\Theta_{c,D})_m,g\rangle$. Summing over $c$ we get $\big\langle(\sum_c\Theta_{c,D})_m,g\big\rangle=0$ for every $D$ and every $x\in\Om$ (the coefficients are rational functions of $x$ and the identity holds for rational functions on $\Om$, hence at every point of $\Om$).
\end{proof}

\emph{The space $I_m$.} Let $I=I^{(n)}\subset\C[V]$ be the ideal of polynomials vanishing on $\bigcup_i(b_iV\cup V_{[n]\setminus i})$. By \eqref{b:ann} and $b_i\Sym^mV=\Sym^m(b_iV)$ we have $I_m=R_m^\perp$ for the above pairing (note $\langle\Psi,u^m\rangle=m!\,\Psi(u)$), that is, $R_m=(I_m)^\perp$. Concretely, for $m\ge2$, by Proposition~\ref{b:dual} and \eqref{b:ann}, $I_m$ is the image in $\C[V]=\C[\nu]/(\sum\nu_C)$ of the span of the monomials $\nu^\chi$ over the \emph{essential} multisets $\chi$ of subsets of $[n]$: $|\chi|=m$ and, for every $i\in [n]$, the number $m_i(\chi)=\#\{C\in\chi:\ i\in C\}$ satisfies $1\le m_i(\chi)\le m-1$ (every $i$ belongs to some, but not to all, sets of $\chi$).
\begin{Proposition}[criterion]\label{inj:crit}
If $I_m\subseteq Z_m(x)$ for some $x\in\Om$, then $(\ast_m)$ holds.
\end{Proposition}
\begin{proof}
If $\widehat G_g\in\Loc^\Om$, then by Proposition~\ref{inj:C} $g\perp Z_m(x)\supseteq I_m$, so $g\in I_m^\perp=R_m$.
\end{proof}

\subsection{Product points}\label{inj:s5}
By Proposition~\ref{inj:crit} it remains to find one point $x$ with $I_m\subseteq Z_m(x)$. For $m=1$ we have $I_1=0$ (the subspaces $b_iV$, $V_{[n]\setminus i}$ span $V$: $\lambda_T\in V_{[n]\setminus i}$ for $i\notin T$, and $\lambda_{[n]}=b_i\lambda_{[n]\setminus i}+\lambda_i$), so $(\ast_1)$ is trivial. For $m\ge2$ we have to show that the monomials $\nu^\chi$ with $\chi$ essential, $|\chi|=m$, lie in $Z_m(x)$ for one point $x$.

\paragraph{The point and the tensor notation.} Fix
\[
x_c=\prod_ip_i(c_i),
\]
where each $p_i$ is a strictly positive probability vector on $\cX_i$ having two letters of different probabilities (for a binary alphabet this means $p_i\ne(\frac12,\frac12)$). Then $N=1$ and $q_{c,C}=\prod_{i\notin C}p_i(c_i)$. Let $U=\bigotimes_i\C^{\cX_i}\cong\C^{\cX}$, with the basis $e_a$ in each factor and $\mathbb 1=\sum_ae_a$. For a cell $c$ and $C\subseteq [n]$ put
\[
\begin{aligned}
1_{c,C}&=\bigotimes_i\big(\mathbb 1\text{ if }i\in C,\ e_{c_i}\text{ if }i\notin C\big)\in U,\\
\ell_{c,C}(u)&=\langle1_{c,C},u\rangle,\qquad
w^c(\nu)=\sum_C\frac{\nu_C}{q_{c,C}}\,1_{c,C}\in U .
\end{aligned}
\]
(Here the $\nu_C$ are treated as \emph{independent} variables; the relation $\sum_C\nu_C=0$ is imposed at the end. Note that $1_{c,\varnothing}=\delta_c$, the basis vector of the cell $c$.) For a homogeneous polynomial $D$ on $U$ of degree $d$ put
\[
T_D(\nu)=\sum_cx_c\,D\big(w^c(\nu)\big),\qquad \Dop=\partial_{\nu_\varnothing}-\tfrac12\sum_C\nu_C\,\partial^2_{\nu_C}
\]
(here $\Dop$ is a second order operator on $\C[\nu_C:\ C\subseteq [n]]$).
\emph{Strategy.} By Lemma~\ref{inj:taylor} below, for $D\in\mathfrak S_k$ the component of degree $k$ of $\sum_c\Theta_{c,D}$ is $T_D$, and for $D\in\mathfrak S_{k+1}$ it is $\Dop T_D$. We choose $D$ of tensor form, for which $T_D$ can be computed (Lemma~\ref{inj:tensor}); the $T_D$ turn out to span the monomials $\nu^\chi$ with $1\le m_i(\chi)\le k-2$ (Lemma~\ref{inj:profile}, Corollary~\ref{inj:mono}); applying $\Dop$ gives the remaining essential monomials, those with some $m_i(\chi)=k-1$ (Lemma~\ref{inj:FC}).

\begin{Lemma}\label{inj:taylor}
(a) For $D\in\mathfrak S_k$ the degree $k$ part (in $\nu$) of $\sum_c\Theta_{c,D}$ is $T_D$. (b) For $D\in\mathfrak S_{k+1}$ it is $\Dop T_D$. (Computed with the $\nu_C$ as independent variables; restrict to $\sum\nu_C=0$ afterwards.)
\end{Lemma}
\begin{proof}
Fix $c$ and put $z_C=\ell_{c,C}(u)/q_{c,C}$. Then $\prod_C(1+z_C)^{\nu_C}=e^{\Lambda_c}$ with
\[
\Lambda_c=\sum_C\nu_C\log(1+z_C)=L-\tfrac12Q+O(u^3),\qquad L=\sum_C\nu_Cz_C=\langle w^c(\nu),u\rangle,\quad Q=\sum_C\nu_Cz_C^2 .
\]
We expand $(x_c+u_c)e^{\Lambda_c}$ by the degree in $u$ and in $\nu$; $L$ has $\nu$-degree $1$ and $u$-degree $1$, $Q$ has $\nu$-degree $1$ and $u$-degree $2$. The two pieces of $\nu$-degree $k$ that can be hit by an operator of degree $k$, resp.\ $k+1$, are
\begin{align*}
&(u\text{-degree }k):\quad x_c\,\frac{L^k}{k!},\\
&(u\text{-degree }k+1):\quad u_c\,\frac{L^k}{k!}-\frac{x_c}2\,\frac{QL^{k-1}}{(k-1)!}.
\end{align*}
Recall: if $D$ is homogeneous of degree $d$ and $l_j(u)=\langle\hat l_j,u\rangle$, then $D(\partial_u)[l_1\cdots l_d]=\partial_{\hat l_1}\cdots\partial_{\hat l_d}D$ (a constant; $\partial_{\hat l}$ is the derivative in the direction $\hat l$). Hence $D(\partial)[l^d]=d!\,D(\hat l)$, $D(\partial)[l_1^{d-1}l_2]=(d-1)!\,(\partial_{\hat l_2}D)(\hat l_1)$ and $D(\partial)[l_1^{d-2}l_2l_3]=(d-2)!\,(\partial_{\hat l_2}\partial_{\hat l_3}D)(\hat l_1)$.

(a) For $D\in\mathfrak S_k$ the first piece gives $x_cD(w^c(\nu))$; summing over $c$ we get $T_D$.

(b) For $D\in\mathfrak S_{k+1}$ the second piece gives
\[
\partial_{\delta_c}D(w^c)-\frac{x_c}2\sum_C\frac{\nu_C}{q_{c,C}^2}\,\partial^2_{1_{c,C}}D(w^c).
\]
Since $w^c(\nu)$ depends linearly on $\nu_C$ with the coefficient $q_{c,C}^{-1}1_{c,C}$, the chain rule gives $\partial_{\nu_C}[D(w^c)]=q_{c,C}^{-1}(\partial_{1_{c,C}}D)(w^c)$ and $\partial^2_{\nu_C}[D(w^c)]=q_{c,C}^{-2}(\partial^2_{1_{c,C}}D)(w^c)$. As $1_{c,\varnothing}=\delta_c$ and $q_{c,\varnothing}=x_c$, we have $\partial_{\delta_c}D(w^c)=x_c\,\partial_{\nu_\varnothing}[D(w^c)]$. So the expression equals $x_c\big(\partial_{\nu_\varnothing}-\frac12\sum_C\nu_C\partial^2_{\nu_C}\big)D(w^c)$, and summing over $c$ gives $\Dop T_D$.
\end{proof}

\paragraph{Operators of tensor form.} Let $K_d(s)$ be the space of linear functionals on $(\C^s)^{\otimes d}$ vanishing on $e_a^{\otimes d}$ ($a=1,\dots,s$) and on $\mathbb 1^{\otimes d}$. For $\varphi_i\in K_d(|\cX_i|)$ put $D(\xi)=(\varphi_1\otimes\dots\otimes\varphi_n)(\xi^{\otimes d})$, regrouping the tensor factors: $U^{\otimes d}=\bigotimes_i(\C^{\cX_i})^{\otimes d}$. Then $D\in\mathfrak S_d$. Indeed, every vector of $W_{i,a}$, resp.\ $W_i$, has the form $\xi=u'\otimes\lambda$ with $\lambda=e_a$, resp.\ $\lambda=\mathbb 1$, in the slot $i$ and an arbitrary vector $u'$ in the other slots; after regrouping, $\xi^{\otimes d}=u'^{\otimes d}\otimes\lambda^{\otimes d}$, so $D(\xi)$ contains the factor $\varphi_i(\lambda^{\otimes d})=0$.
\begin{Lemma}\label{inj:tensor}
For such $D$, identify a $d$-tuple of subsets $(C_1,\dots,C_d)$ with $(A_1,\dots,A_n)$, $A_i=\{t\in[d]:\ i\in C_t\}$. Then
\[
T_D(\nu)=\sum_{(A_i)}\nu_{C_1}\cdots\nu_{C_d}\prod_iH^{(i)}(A_i),\qquad H^{(i)}(A)=\sum_{a\in \cX_i}p_i(a)^{\,1-(d-|A|)}\varphi_i\big(e_a^{\otimes A^c}\otimes\mathbb 1^{\otimes A}\big),
\]
where $e_a^{\otimes A^c}\otimes\mathbb 1^{\otimes A}$ denotes the tensor with $e_a$ in the positions $t\notin A$ and $\mathbb 1$ in the positions $t\in A$.
\end{Lemma}
\begin{proof}
Expand $w^c(\nu)^{\otimes d}$ over the tuples $(C_1,\dots,C_d)$. The weight of a tuple is $\prod_t\nu_{C_t}/q_{c,C_t}=\prod_t\nu_{C_t}\cdot\prod_ip_i(c_i)^{-(d-|A_i|)}$, since $q_{c,C_t}=\prod_{i\notin C_t}p_i(c_i)$ and $i\notin C_t$ for exactly $d-|A_i|$ values of $t$. After regrouping the factors, the tensor $1_{c,C_1}\otimes\dots\otimes1_{c,C_d}$ becomes $\bigotimes_i\big(e_{c_i}^{\otimes A_i^c}\otimes\mathbb 1^{\otimes A_i}\big)$, so $\bigotimes_i\varphi_i$ takes the value $\prod_i\varphi_i\big(e_{c_i}^{\otimes A_i^c}\otimes\mathbb 1^{\otimes A_i}\big)$. Multiply by $x_c=\prod_ip_i(c_i)$ and sum over $c$; the sum factorizes over $i$.
\end{proof}
\begin{Lemma}\label{inj:profile}
Let $p$ be a strictly positive probability vector on a finite set $\cX_0$ with $p_a\ne p_b$ for some $a,b$, and let $d\ge2$. As $\varphi$ ranges over $K_d(|\cX_0|)$, the functions $H_\varphi(A)=\sum_ep_e^{1-(d-|A|)}\varphi(e_e^{\otimes A^c}\otimes\mathbb 1^{\otimes A})$ on the subsets $A\subseteq[d]$ range over \emph{all} functions supported on $\{A:\ 1\le|A|\le d-2\}$.
\end{Lemma}
\begin{proof}
\emph{Idea.} $H_\varphi$ is linear in $\varphi$; the two constraints $\varphi(e_e^{\otimes d})=0$, $\varphi(\mathbb 1^{\otimes d})=0$ force $H_\varphi$ to vanish for $|A|\in\{0,d-1,d\}$. We show that there are no other constraints.

\emph{(1) The support.} $H(\varnothing)=\sum_ep_e^{1-d}\varphi(e_e^{\otimes d})=0$; for $|A|=d-1$,
\[
\sum_e\varphi(e_e\text{ at one place},\ \mathbb 1\text{ elsewhere})=\varphi(\mathbb 1^{\otimes d})=0
\]
because $\sum_ee_e=\mathbb 1$ (here the weight is $p_e^{0}=1$); and $H([d])=\sum_ep_e\varphi(\mathbb 1^{\otimes d})=0$. So $H_\varphi$ is supported on $1\le|A|\le d-2$.

\emph{(2) Reduction to two letters.} Let $\pi:\C^{\cX_0}\to\C^2$, $e_a\mapsto f_0$, $e_b\mapsto f_1$, $e_e\mapsto0$ for the other letters, and $\varphi=\psi\circ\pi^{\otimes d}$ with $\psi\in K_d(2)$ (this $\varphi$ lies in $K_d(|\cX_0|)$: $\pi(\mathbb 1)=f_0+f_1=:\mathbb 1_2$, so $\varphi$ vanishes on $e_e^{\otimes d}$ and on $\mathbb 1^{\otimes d}$ because $\psi$ vanishes on $f_0^{\otimes d}$, $f_1^{\otimes d}$, $\mathbb 1_2^{\otimes d}$). For $A^c\ne\varnothing$ the terms with $e\ne a,b$ vanish, so with $p=p_a$, $q=p_b$, $s=d-|A|$:
\[
H_\varphi(A)=p^{1-s}\psi(f_0^{\otimes A^c}\otimes\mathbb 1_2^{\otimes A})+q^{1-s}\psi(f_1^{\otimes A^c}\otimes\mathbb 1_2^{\otimes A}).
\]
It remains to show that these two-letter functions fill the space of functions supported on $1\le|A|\le d-2$; only $p\ne q$ (both positive) is used, no relation such as $p+q=1$.

\emph{(3) Two letters: coordinates.} The tensors $w_T=\bigotimes_t(\mathbb 1_2\text{ if }t\in T,\ f_0\text{ if }t\notin T)$, $T\subseteq[d]$, form a basis of $(\C^2)^{\otimes d}$; put $\varphi_T=\psi(w_T)$ (these $2^d$ numbers are arbitrary coordinates of $\psi$). Then $f_0^{\otimes d}=w_\varnothing$, $\mathbb 1_2^{\otimes d}=w_{[d]}$, $f_0^{\otimes A^c}\otimes\mathbb 1_2^{\otimes A}=w_A$, and, because $f_1=\mathbb 1_2-f_0$ and we expand, $f_1^{\otimes A^c}\otimes\mathbb 1_2^{\otimes A}=\sum_{T\supseteq A}(-1)^{d-|T|}w_T$. Put
\[
\Psi(A)=\psi(f_1^{\otimes A^c}\otimes\mathbb 1_2^{\otimes A})=\sum_{T\supseteq A}(-1)^{d-|T|}\varphi_T,\qquad\text{so that}\qquad\varphi_A=\sum_{T\supseteq A}(-1)^{d-|T|}\Psi(T)
\]
(M\"obius inversion on the Boolean lattice of subsets of $[d]$). Then $\psi\in K_d(2)$ means $\varphi_\varnothing=0$, $\Psi(\varnothing)=\psi(f_1^{\otimes d})=0$, $\Psi([d])=\varphi_{[d]}=0$. Since $\varphi_\varnothing=\sum_T(-1)^{d-|T|}\Psi(T)=:R(\Psi)$, the constraints are: $\Psi(\varnothing)=\Psi([d])=0$ and $R(\Psi)=0$. In terms of $\Psi$,
\[
H(A)=q^{1-s}\Psi(A)+p^{1-s}\sum_{T\supseteq A}(-1)^{d-|T|}\Psi(T),\qquad s=d-|A| .
\]
\emph{(4) Solving the triangular system.} The term $T=A$ of the sum is $(-1)^s\Psi(A)$, so $H(A)=c_s\Psi(A)+p^{1-s}\sum_{T\supsetneq A}(-1)^{d-|T|}\Psi(T)$ with $c_s=q^{1-s}+(-1)^sp^{1-s}$. For $s=1$ (i.e.\ $|A|=d-1$) we get $c_1=0$ and $H(A)=0$ automatically (using $\Psi([d])=0$), in accordance with (1). For $s\ge2$ we have $c_s\ne0$: it is positive for even $s$, and for odd $s$ it vanishes iff $p=q$. Therefore, given a target $h$ supported on $1\le|A|\le d-2$ and arbitrary values $a_j=\Psi([d]\setminus j)$, $j\in[d]$, the equations $H(A)=h(A)$ for $|A|=d-2,d-3,\dots,1$ determine $\Psi(A)$ step by step (the equation for $A$ involves only $\Psi(A)$ and the values at strictly larger sets), and $\Psi(\varnothing)=0$, $\Psi([d])=0$ are imposed separately; the equations with $|A|\ge1$ do not involve $\Psi(\varnothing)$.

\emph{(5) The last constraint $R(\Psi)=0$.} The solution is linear in $(h,a)$, so $R(\Psi)=\rho(h)+\sum_j\kappa_ja_j$ with a linear form $\rho$. The whole problem is invariant under permutations of $[d]$, which permute the $a_j$; hence all $\kappa_j$ are equal to one number $\kappa$. If $\kappa\ne0$, we choose $a$ with $\kappa\sum_ja_j=-\rho(h)$ and get a solution with $R(\Psi)=0$, that is $\psi\in K_d(2)$ with $H_\varphi=h$. To compute $\kappa$ take $h=0$ and all $a_j=1$ (then $d\kappa=R(\Psi)$). The solution of the homogeneous system with these $a_j$ is
\[
\Psi(T)=u_{d-|T|}\quad(T\ne\varnothing),\qquad u_j=j\,r^{j-1},\quad r=\frac q{p+q},\qquad\Psi(\varnothing)=0 .
\]
Indeed, for $1\le|A|\le d-2$ we have $\sum_{T\supseteq A}(-1)^{d-|T|}\Psi(T)=\sum_{j=0}^{s}\binom sj(-1)^ju_j=-s(1-r)^{s-1}$ (here $u_0=0$ and we used $\sum_j\binom sj(-1)^jjr^{j-1}=-s(1-r)^{s-1}$), so
\[
H(A)=q^{1-s}sr^{s-1}-p^{1-s}s(1-r)^{s-1}=s(p+q)^{1-s}-s(p+q)^{1-s}=0 .
\]
Finally $R(\Psi)=\sum_{s=1}^{d-1}\binom ds(-1)^ssr^{s-1}$ (the terms $s=0$ and $s=d$ are absent since $\Psi([d])=\Psi(\varnothing)=0$), and the same binomial identity gives
\[
R(\Psi)=-d\Big[(1-r)^{d-1}-(-r)^{d-1}\Big]\ne0\qquad(p\ne q):
\]
for $d-1$ odd both terms in the bracket are positive, for $d-1$ even it equals $(1-r)^{d-1}-r^{d-1}\ne0$ because $r\ne\frac12$. Hence $\kappa\ne0$.
\end{proof}
\begin{Corollary}\label{inj:mono}
The polynomials $T_D$, $D=\bigotimes\varphi_i$ with $\varphi_i\in K_d(|\cX_i|)$, span $M_d:=\operatorname{span}\{\nu^\chi:\ |\chi|=d,\ 1\le m_i(\chi)\le d-2\ \forall i\}$ (formal monomials in the $2^n$ variables $\nu_C$).
\end{Corollary}
\begin{proof}
By Lemma~\ref{inj:tensor} the map $(\varphi_i)\mapsto T_D$ is the image of $\bigotimes_iH^{(i)}$ under $h\mapsto\sum_{(A_i)}h(A_1,\dots,A_n)\,\nu^{\chi(A)}$, where $\chi(A)=\{C_1,\dots,C_d\}$. By Lemma~\ref{inj:profile} each $H^{(i)}$ can be the indicator function of any single set $A_i$ with $1\le|A_i|\le d-2$; this gives the monomial $\nu^{\chi(A)}$ with $m_i(\chi(A))=|A_i|$. Conversely, all $T_D$ lie in $M_d$ by the support statement of Lemma~\ref{inj:profile}.
\end{proof}
\begin{Lemma}\label{inj:FC}
Let $M_k^{+}=M_k+\Dop M_{k+1}\subseteq\C[\nu_C]_k$. Every $\nu^\chi$ with $\chi$ essential, $|\chi|=k\ge2$, lies in $M_k^{+}$.
\end{Lemma}
\begin{proof}
For a monomial $\nu^\omega$ with $a=m_\varnothing(\omega)$ and $m_C=m_C(\omega)$ (the multiplicity of $C$), $\partial_{\nu_\varnothing}\nu^\omega=a\,\nu^{\omega-\varnothing}$ and $-\frac12\nu_C\partial^2_{\nu_C}\nu^\omega=-\frac12m_C(m_C-1)\nu^{\omega-C}$ (removal of one copy). Summing over $C$ (including $C=\varnothing$, which contributes $-\frac12a(a-1)\nu^{\omega-\varnothing}$),
\[
\Dop\nu^\omega=\frac{a(3-a)}2\,\nu^{\omega-\varnothing}-\frac12\sum_{C\ne\varnothing}m_C(m_C-1)\,\nu^{\omega-C}.
\]
We use induction on the weight $w(\chi)=\sum_{C\in\chi}|C|$ (an essential $\chi$ has $w\ge n\ge1$). If $m_i(\chi)\le k-2$ for all $i$, then $\nu^\chi\in M_k$. Otherwise $m_i(\chi)=k-1$ for some $i$: exactly one element of $\chi$ avoids $i$, and $\varnothing$ avoids $i$, so $a_0:=m_\varnothing(\chi)\le1$. Put $\omega=\chi+\varnothing$ (adding one copy of $\varnothing$; this does not change the weight). Then $m_j(\omega)=m_j(\chi)\in[1,(k+1)-2]$ for all $j$ (here $m_j(\chi)\le k-1$ because $\chi$ is essential), so $\nu^\omega\in M_{k+1}$. Moreover $m_\varnothing(\omega)=a_0+1\in\{1,2\}$, so the coefficient of $\nu^\chi$ in $\Dop\nu^\omega$ is $\frac{(a_0+1)(2-a_0)}2=1$. The other terms are $-\frac12m_C(m_C-1)\nu^{\chi'}$ with $\chi'=\omega-C$, where $C\ne\varnothing$ occurs at least twice in $\omega$ (otherwise the coefficient is $0$), hence at least twice in $\chi$. Such a $\chi'$ is essential ($C$ occurs twice in $\chi$, so $m_j(\chi')\ge1$; and $m_j(\chi')\le m_j(\chi)\le k-1$) and has $w(\chi')<w(\chi)$. By induction $\nu^{\chi'}\in M_k^{+}$, hence $\nu^\chi=\Dop\nu^\omega+\frac12\sum m_C(m_C-1)\nu^{\chi'}\in M_k^{+}$.
\end{proof}
\begin{Proposition}[covering at product points]\label{inj:ZM}
For $x$ as above and $k\ge2$, the images in $\C[V]$ of all $\nu^\chi$ with $\chi$ essential, $|\chi|=k$, lie in $Z_k(x)$. Hence $I_m\subseteq Z_m(x)$ for all $m$.
\end{Proposition}
\begin{proof}
For $D\in\mathfrak S_k$, resp.\ $\mathfrak S_{k+1}$, of the form $\bigotimes_i\varphi_i$, the degree $k$ component of $\sum_c\Theta_{c,D}$ is $T_D$, resp.\ $\Dop T_D$ (Lemma~\ref{inj:taylor}); so $Z_k(x)\supseteq M_k+\Dop M_{k+1}$ (Corollary~\ref{inj:mono}) after restriction to $\sum\nu_C=0$. By Lemma~\ref{inj:FC} this contains all essential $\nu^\chi$. By Proposition~\ref{b:dual} these span $I_k$ for $k\ge2$, and $I_1=0$.
\end{proof}

\subsection{Conclusion}
\begin{proof}[Proof of Proposition~\ref{inj:A}]
Choose $x$ as in \S\ref{inj:s5}; this is possible for arbitrary alphabets $|\cX_i|\ge2$ (take every $p_i$ positive and non-uniform). Proposition~\ref{inj:ZM} gives $I_m\subseteq Z_m(x)$, Proposition~\ref{inj:crit} gives $(\ast_m)$, and by \eqref{inj:ker} and Lemma~\ref{inj:red} the kernel of $Q_n^{(m)}\to H^n(S,\FF)$ is contained in $\{g:\widehat G_g\in\Loc^\Om\}/R_m$, which is $0$ by $(\ast_m)$.
\end{proof}
\begin{Remark}
(1) The argument as written needs two letters of different probability in every $p_i$. For $p_i=(\frac12,\frac12)$ one has $c_s=0$ for all odd $s$, so that for $d\ge4$ the triangular system in step (4) of the proof of Lemma~\ref{inj:profile} cannot be solved as described (for $d=3$ it can, but $\kappa=0$), and the conclusion of that lemma fails for every $d\ge3$. We do not claim that the inclusion $I_m\subseteq Z_m(x)$ fails at such points: $Z_m(x)$ involves operators of all degrees, whereas the argument uses only the degrees $m$ and $m+1$. For $|\cX_i|\ge3$ a uniform $p_i$ is not covered by the proof. Since the criterion needs just one point $x$, no generality is lost.
(2) Of the space $\mathfrak S_d$ we only use that it contains the operators built from the $K_d$ (not a description of all of $\mathfrak S_d$); no regularity of the functions in $\FF$ is assumed (Proposition~\ref{inj:kill} uses finite differences).
(3) Only the constant term in $z$ (Proposition~\ref{inj:C}) and operators of degrees $k$ and $k+1$ for the target degree $k=m$ are used.
\end{Remark}

\section{All degrees; degree one}\label{sec:cons}

\begin{Proposition}[the entropy class]\label{inj:H1}
Let $n\ge2$ and $|\cX_i|\ge2$. (a) Every coboundary in degree $1$ of the reduced complex $\tilde C^\bullet(\{1\},\dots,\{n\};\FF)$ is zero; hence $H^1(S,\FF)$ is the space of cocycles, that is, of families $(f_i)_{i\in [n]}$, $f_i\in\FF(\{i\})$, $f_i(\delta_a)=0$ for point masses $\delta_a$, such that $f_i(P_i)-\sum_bP_j(b)\,f_i\big((P|_{X_j=b})_i\big)=f_j(P_j)-\sum_aP_i(a)\,f_j\big((P|_{X_i=a})_j\big)$ for all $i\ne j$. (b) The generator of $H^1(S,\LL^{(1)})=\R$ is mapped by $\ev_*$ to the class of the Shannon entropies $(H(P_i))_i$, which is nonzero. Consequently $\ev_*:H^1(S,\LL^{(1)})\to H^1(S,\FF)$ is injective, and its image is the line spanned by the entropy class.
\end{Proposition}
\begin{proof}
We use Proposition~\ref{cor:reduced} for the atoms $\{1\},\dots,\{n\}$: $\tilde C^0=\FF(\varnothing)$, $\tilde C^1=\bigoplus_ie_i\FF(\{i\})$, and the differential sends $u\in\tilde C^0$ to $(e_i\iota u)_i$ and $(f_i)$ to $\big(e_ie_j(\iota f_j-\iota f_i)\big)_{i<j}$ up to sign. (a) $\FF(\varnothing)=\R$ consists of constants (there is a single law on the empty product), and for a constant $u$ we have $b_iu=\sum_aP_i(a)u=u$, so $e_i\iota u=0$. Since $e_i\iota f_i=\iota e_if_i=\iota f_i$, the cocycle condition reads $e_j\iota f_i=e_i\iota f_j$, which is the displayed equation; and $f_i\in e_i\FF(\{i\})$ means $f_i(\delta_a)=0$ for all $a$, because $(b_if)(P)=\sum_aP(a)f(\delta_a)$ on $\FF(\{i\})$.
(b) Since $\LL^{(1)}(\varnothing)=0$, there are no coboundaries in $\LL^{(1)}$ either. The family $f_i=\lambda_i$ is a cocycle ($b_i\lambda_i=0$, so $f_i=e_i\lambda_i$; and $e_j\lambda_i=\lambda_i+\lambda_j-\lambda_{ij}=e_i\lambda_j$), it is nonzero, hence it spans $H^1(S,\LL^{(1)})$, which is one-dimensional by Theorem~\ref{thm:conc}. Its image under $\ev$ is $(H(P_i))_i$, a nonzero family of functions; by (a) it is not a coboundary.
\end{proof}
\begin{Corollary}[injectivity in all degrees]\label{inj:all}
For $n\ge2$, arbitrary finite alphabets $|\cX_i|\ge2$ and all $k\ge0$, $m\ge1$, the map $\ev_*:H^k(S,\LL^{(m)})\to H^k(S,\FF)$ is injective.
\end{Corollary}
\begin{proof}
By the concentration theorem (Theorem~\ref{thm:conc}) $H^k(S,\LL^{(m)})=0$ unless $k=n$, or $k=1$ and $m=1$. The case $k=n$ is Proposition~\ref{inj:A}, the case $(k,m)=(1,1)$ is Proposition~\ref{inj:H1}.
\end{proof}
\begin{proof}[Proof of Theorem~\ref{thm:A}]
By Lemma~\ref{inj:Pi}(c) the map $\ev_*:H^k(S,\LL^{(m)})\to H^k(S,\FF)$ is the composite of $(\ev_\Pi)_*$ and $\mathrm{res}_*$. It is injective by Corollary~\ref{inj:all}, hence so is its first factor.
\end{proof}
\subsection{Degree one with measurable coefficients}

\begin{Proposition}\label{inj:H1V}
For $n\ge2$ and $|\cX_i|\ge2$ the map $(\ev_\Pi)_*:H^1(S,\LL^{(1)})\to H^1(S,\FF_\Pi)$ is an isomorphism onto the line spanned by the Shannon entropy.
\end{Proposition}
\begin{proof}
It is injective by Theorem~\ref{thm:A}, and $\dim H^1(S,\FF_\Pi)=1$ by Vigneaux's theorem \cite[Theorem~1.3.2]{V}, whose hypotheses hold: all laws are allowed, so the marginalizations are surjective; the objects with at least two outcomes are the nonempty $y\subseteq [n]$ and form a connected category; and each of them is refined by $[n]=\{i\}\cup([n]\setminus i)$, a non-degenerate product in the sense of \cite[Def.~4.5.6]{V}. The proof of that theorem rests on the description of the measurable solutions of the fundamental equation of information \cite{Tv,KN}. The statement is also part~(ii) of \cite[Proposition~16]{LKS1} and is repeated here for completeness; the only addition to Vigneaux's theorem is that the formal class of $(\lambda_i)$ is mapped to a generator.
\end{proof}

\begin{Remark}[arbitrary functions give a larger group]\label{inj:ell}
For a homomorphism $\ell:(\R_{>0},\cdot)\to(\R,+)$ put $H_\ell(P)=-\sum_aP(a)\,\ell(P(a))$ (sum over the support). Since $\ell(ab)=\ell(a)+\ell(b)$, the chain rule $H_\ell(P_{ij})=H_\ell(P_i)+\sum_aP_i(a)H_\ell(P_{j\mid a})$ holds, and $H_\ell$ vanishes on point masses; hence $(H_\ell(P_i))_i$ is a cocycle in the sense of Proposition~\ref{inj:H1}(a) and defines a class in $H^1(S,\FF)$. The logarithms of $2$ and $3$ are linearly independent over $\mathbb Q$, so one can choose $\ell$ with $\ell(2)=1$ and $\ell(3)=0$. Then $H_\ell(\frac12,\frac12)=1$ while $H_\ell(\frac13,\frac23)=-\frac23$. The Shannon entropy is positive on both laws, so $H_\ell$ is not a multiple of it; this already happens for binary alphabets. Hence $H^1(S,\FF)$ is strictly larger than the image of $H^1(S,\LL^{(1)})$. Such $\ell$ are not measurable, and these classes do not occur in $H^1(S,\FF_\Pi)$.
\end{Remark}

\subsection{Moments}

\begin{proof}[Proof of Corollary~\ref{cor:B}]
(i)$\Rightarrow$(ii): write $g=\sum_ib_ih_i+\sum_ik_i$ with $h_i\in\Sym^mV$ and $k_i\in\Sym^mV_{[n]\setminus i}$. Since $\ev_\Pi$ is a morphism of modules (Lemma~\ref{inj:Pi}), \eqref{eq:dec} holds on all laws with $\varphi_i=\ev_\Pi(h_i)$ and $\psi_i=\ev_\Pi(k_i)$. (ii)$\Rightarrow$(iii) is clear. (iii)$\Rightarrow$(i): both sums in \eqref{eq:dec} involve, for positive $P$, only functions of slices and of marginals of $P$, so the computation in the proof of Lemma~\ref{inj:red} shows that $\widehat G_g\in\Loc^\Om$. Now apply $(\ast_m)$, which is established in the proof of Proposition~\ref{inj:A}.
\end{proof}

For example, for every $m\ge2$ the class of $\lambda_{[n]}^m$ in $Q_n^{(m)}$ is nonzero {\cite[Theorem~6]{LKS1}}. Hence the $m$-th moment $P\mapsto\EE_P\,(-\log P(X_{[n]}))^m$ of the self-information of the whole system admits no decomposition as in Corollary~\ref{cor:B}. For $m=2$ this is the linear part of the varentropy; the term $H(P)^2$ is not a moment and is not covered.

\section{Open questions}\label{sec:open}

\begin{enumerate}[leftmargin=2em]
\item Compute $H^k(S,\FF_\Pi)$ for $k\ge2$. For which classes of regular functions is the sum over $m$ of the images of $H^n(S,\LL^{(m)})$ the whole top cohomology?
\item Is there a proof of Proposition~\ref{inj:A} for measurable functions which avoids finite differences and the choice of a product law?
\item Can products of entropies, and hence varentropy and covariances, be included by a suitable enlargement of the coefficient modules?
\end{enumerate}

\subsection*{Statements} No datasets were generated or analyzed. \textbf{Author contributions and conflicts of interest: to be completed and confirmed by the authors.}

\section*{Acknowledgements}
The work was carried out with the support of the federal assignment of the Ministry of Science and Higher Education of the Russian Federation, project No.~FSMG-2025-0070.

\end{document}